\documentclass[11pt]{article}

\usepackage{graphics}

\usepackage{latexsym}
\usepackage{amsmath}
\usepackage{amssymb}
\usepackage[dvips]{epsfig}

\newtheorem{theorem}{Theorem}[section]

\newtheorem{lemma}{Lemma}[section]

\newcommand{\bt}{\begin{theorem}}
\newcommand{\bl}{\begin{lemma}}
\newcommand{\el}{\end{lemma}}
\newcommand{\et}{\end{theorem}}

\newcommand{\bn}{\begin{eqnarray}}
\newcommand{\en}{\end{eqnarray}}
\newcommand{\bnn}{\begin{eqnarray*}}
\newcommand{\enn}{\end{eqnarray*}}

\makeatletter      % '@' is now a normal "letter" for TeX
\@addtoreset{equation}{section}
\makeatother       % '@' is restored as a "non-letter" character for TeX

\title{On the Inviscid Limit of the 3D Navier-Stokes Equations
with Generalized Navier-slip Boundary Conditions\footnote{This research is
supported in part by NSFC 10971174, and Zheng Ge Ru Foundation,
and Hong Kong RGC Earmarked Research Grants CUHK-4041/11P,
CUHK-4042/08P, a Focus Area Grant from The Chinese University of Hong Kong, and a grant from Croucher Foundation.}}

\author{Yuelong Xiao\footnote{Institute for Computational
and Applied mathematics, Xiangtan University.}$\ ^{, \ddag}$ and
Zhouping Xin\footnote{The Institute of Mathematical Sciences, The
Chinese University of Hong Kong.}}

\date{}

\begin{document}
\maketitle

\vskip 1cm

\begin{abstract}
In this paper, we investigate the vanishing viscosity limit problem for the 3-dimensional (3D)
incompressible Navier-Stokes equations in a general bounded smooth domain of $R^3$
with the generalized Navier-slip boundary conditions (\ref{VSg}).
Some uniform estimates on rates of convergence in $C([0,T],L^2(\Omega))$ and
$C([0,T],H^1(\Omega))$ of the solutions to the corresponding solutions
of the idea Euler equations with the standard slip boundary condition are obtained.
\end{abstract}

\vskip 1cm

\section{Introduction}

Let $\Omega\subset R^3$ be a bounded smooth domain. We consider the Navier-Stokes
equations
\begin{eqnarray}\label{NS1}
&& \partial_{t}u^\varepsilon -\varepsilon \Delta u^\varepsilon
   + u^\varepsilon\cdot\nabla u^\varepsilon  + \nabla p^\varepsilon = 0,
\textrm{ in }\Omega\\ \label{NS2}
&& \nabla\cdot u^\varepsilon  = 0, \textrm{ in }\Omega
\end{eqnarray}
and the Euler equations
\begin{eqnarray}\label{Eu}
&& \partial_{t}u + u\cdot\nabla u  + \nabla p = 0,
\textrm{ in }\Omega\\\label{Eu1}
&& \nabla\cdot u  = 0, \textrm{ in }\Omega
\end{eqnarray}

We are interested in the vanishing viscosity limit problem:
do the solutions of the Navier-Stokes equations (\ref{NS1}),(\ref{NS2})
converge to that of the Euler equations (\ref{Eu}) with the same initial data
as $\varepsilon$ vanishes?.

The vanishing viscosity limit problem for the Navier-Stokes equations is a
classical issue and has been well studied when the domain has no boundaries, various
convergence results have been obtained, see for instances\cite{Co,Co1,EM,Ka,Ka1,Mas,Sw}.

However, in the presence of a physical boundary, the situations is much more complicated,
and the problem becomes challenging due to the boundary layers. This is so even in
the case that the corresponding solution to the initial boundary value problem of the
Euler equations remains smooth.

Indeed, in presence of a boundary $\partial\Omega$, one of the most important physical
boundary conditions for the Euler equations (\ref{Eu}) is the following slip boundary condition, i.e.,
\begin{equation}\label{Sl}
 u\cdot n = 0\ {\rm on} \ \partial \Omega
\end{equation}
and the initial boundary value problem of the Euler equations (\ref{Eu}) with the slip
boundary condition (\ref{Sl}) has a smooth solution at least locally in time. Corresponding
to the slip boundary condition (\ref{Sl}) for the Euler equations, there are
different choices of boundary conditions for the Navier-Stokes equations,
one is the mostly-used no-slip boundary condition, i.e.,
\begin{equation}\label{nSl}
  u^\varepsilon  = 0,\ {\rm on}\ \Omega
\end{equation}
This Dirichlet type boundary condition was proposed by Stokes (\cite{St}) assuming
that fluid particles are adherent to the boundary due to the positive viscosity.
Although the well-posedness of the smooth solution to the initial boundary value
problem of the Navier-Stokes equations with the no-slip boundary condition can be established
quiet easily (at least local in time), the asymptotic convergence of the solution
to the corresponding solution of the Euler equations (\ref{Eu}) with the boundary condition
(\ref{Sl}) as the viscosity coefficient $\varepsilon$ tends to zero is one of major
open problems except special cases (see \cite{Co2,Maz,Sm1,Sm2}) due to the possible appearance of the boundary
layers. In general, only some sufficient conditions
are obtained for the convergence in $L^2(\Omega)$, see \cite{Ka2,Ke,Wa}. However, there are
some new results in 2D case announced recently in \cite{YM}.

Another class of boundary conditions for the Navier-Stokes equations is the Navier-slip boundary conditions, i.e,
\begin{equation}\label{NaS}
   u^\varepsilon\cdot n = 0, \ [2(S(u^\varepsilon)n) + \gamma u^\varepsilon]_{\tau} = 0, \ {\rm on} \ \partial\Omega
    \end{equation}
where $ 2S(u^\varepsilon) = (\nabla u^\varepsilon + (\nabla u^\varepsilon)^T)$ is the viscous stress tensor, and
$\gamma$ is a smooth given function. Such conditions were proposed by Navier in \cite{Na}, which allow the fluid to slip at the boundary,
and have important applications for problems with rough boundaries,
prorated boundaries, and interfacial boundary problems, see \cite{B67,Ber,PT,TT}.
The well-posedness of the Naver-Stokes equations (\ref{NS1}), (\ref{NS2})
with the Navier-slip boundary conditions and related results have been established,
see \cite{Ber,XX1} and the references therein.

In a recent paper \cite{Ke2}, the Navier-slip boundary condition (\ref{NaS}) is written to
the following generalized one
\begin{equation}\label{Nag}
   u^\varepsilon\cdot n = 0, \ [2(S(u^\varepsilon)n) + A u^\varepsilon]_{\tau} = 0, \ {\rm on} \ \partial\Omega
    \end{equation}
with $ A$ is a smooth symmetric tensor. The well-posedness of the Naver-Stokes equations (\ref{NS1}), (\ref{NS2})
with the generalized Navier-slip boundary condition is indeed similar to the classical one.

On the other hand, the following vorticity-slip boundary condition
\begin{equation}\label{VS}
 u^\varepsilon\cdot n = 0,\  n\times(\omega^\varepsilon) = \beta u^\varepsilon \ {\rm on} \
\partial \Omega
\end{equation}
with $\beta$ is a smooth function is introduced in \cite{XX1}, where $\omega = \nabla\times u$ is the vorticity
of the fluid.

It is noticed that
\begin{equation}
   (2(S(v)n) - (\nabla\times v)\times n)_\tau = GD(v)_\tau \ {\rm on} \ \partial\Omega
    \end{equation}
where $GD(v) = -2S(n)v$ (see\cite{XX1}). Hence the generalized slip condition (\ref{Nag})
is equivalent to the following generalized vorticity-slip condition
\begin{equation}\label{VSg}
 u^\varepsilon\cdot n = 0,\  n\times(\omega^\varepsilon) = [B u^\varepsilon]_\tau \ {\rm on} \
\partial \Omega
\end{equation}
with $B$ a given smooth symmetric tensor on the boundary.
The equivalence can also be seen for instance in \cite{BBer,BF,BF1,BN10,WWX,XX}.

Compared to the case of no-slip boundary
condition, the asymptotic behavior of solutions to the Navier-Stokes equations with the
Navier-slip boundary conditions as $\varepsilon\rightarrow 0$ is relatively easy to
establish. Some strong convergence results in 2D (see \cite{CMR,Li}) and various weak
convergence results in 3D (see \cite{IfP,Ke1,Ke2}) are also established, where it has been shown that
\begin{equation}\label{ke1}
  \|u^\varepsilon-u\|^{2} +
  \varepsilon\int_0^t\|u^\varepsilon-u\|_1^{2}dt
 \leq c\varepsilon^{\frac{3}{2}}
 \end{equation}
and
\begin{equation}\label{ke2}
  \|u^\varepsilon-u\|_{L^\infty([0,T]\times\Omega)}
 \leq c\varepsilon^{\frac{3}{8}(1-s)}
 \end{equation}
for some $s>0$, see \cite{Ke2} for details.

For the homogeneous case, i.e, $B=0$ in (\ref{VSg}), better convergence
results are available. The $H^3(\Omega)$ convergence and $H^2(\Omega)$ estimate on
the rate of convergence are obtained in \cite{XX} for flat domains.
These results are improved to $W^{k,p}(\Omega)$ in \cite{BF,BF1}.
However, such a strong convergence can not be expected for general domains,
except the case that the initial vorticity vanishes on the boundary \cite{XX2}.
Note that the convergence in $H^2(\Omega)$ or, $W^{k,p}(\Omega)$ (for $k\geq 2$),
implies that limiting solution of the Euler equation satisfies the boundary
condition (\ref{VSg}) with $B =0$.

However, it is shown in \cite{BF2,BF3}
that the solution to the Euler equations with the slip boundary condition
(\ref{Sl}) can not satisfy the extra condition $n\times\omega=0$ on the
boundary in general, and satisfy $n\times\omega=0\ \rm{on}\ \partial\Omega$
only if $\omega\cdot n = 0\ \rm{on}\ \partial\Omega$. This implies that the only
possible case for convergence of solutions to the Navier-Stokes equations in $H^2(\Omega)$
is the one obtained in \cite{XX2}, where the following estimate on the convergence rate
\begin{equation}\label{xx0}
  \|u^\varepsilon-u\|^{2}+\varepsilon\|u^\varepsilon-u\|_1^{2} +
  \varepsilon\int_0^t(\|u^\varepsilon-u\|_1^{2}+\varepsilon\|u^\varepsilon-u\|_2^{2})dt
 \leq c\varepsilon^{2}
 \end{equation}
is obtained. Furthermore, a $W^{s,p}(\Omega)$ convergence result can
also be obtained in this case, see \cite{Ber1}.

For general initial data and domains in the homogeneous case, i.e, $B=0$ in (\ref{VSg}),
the best estimate on the rate of convergence established so far is in $W^{1,p}(\Omega)$
in \cite{WXZ}, which implies in the particular case $p=2$ that
\begin{equation}\label{wxz}
  \|u^\varepsilon-u\|^{2}+\varepsilon\|u^\varepsilon-u\|_1^{2}
 \leq c\varepsilon^{2-s}
 \end{equation}
with $s =\frac{1}{2}$, where $u^\varepsilon$ and $u$ are solutions to
(\ref{NS1}), (\ref{NS2}), (\ref{VS}) and (\ref{Eu}),
(\ref{Sl}) respectively.

However, all the strong convergence results depend essentially the
homogeneous property, i.e, $n\times\omega = 0$ on $\partial\Omega$.
Indeed, for general $B$(not identically equal to zero on $\partial\Omega$),
there is no any estimate on the rate of convergence in $H^{k}(\Omega)$
(or $W^{k,p}(\Omega)$) for $k\geq1$, as far as the authors are aware.
It should be noted that the co-normal
uniform estimates have been established in \cite{Mas1} for general domains
and general Navier-slip boundary conditions,
which grantee the convergence of solutions to the Navier-Stokes equations
to that of the Euler equations in $W^{1,\infty}(\Omega)$. Yet, some additional
efforts are still needed to obtain an estimate on the rate of convergence
as in (\ref{xx0}).

In a recent paper \cite{XX3}, we proposed the following slip boundary condition
\begin{equation}\label{xxn}
 u^\varepsilon\cdot n = 0,\ \omega^\varepsilon\cdot n = 0, n\times(\Delta u^\varepsilon) = 0\ {\rm on} \
\partial \Omega
\end{equation}
for the Navier-Stokes equations (\ref{NS1}), (\ref{NS2}) and obtained the following estimate on the rate of convergence
\begin{equation}\label{xxnc}
  \|u^\varepsilon-u\|_1^{2}+\varepsilon\|u^\varepsilon-u\|_2^{2} +
  \varepsilon\int_0^t(\|u^\varepsilon-u\|_2^{2}+\varepsilon\|u^\varepsilon-u\|_3^{2})dt
 \leq c\varepsilon^{2-s}
 \end{equation}
for any $s>0$.

In this paper, we investigate the vanishing
viscosity problem for the Navier-Stokes equations (\ref{NS1}), (\ref{NS2}) on
a general bounded smooth domain $\Omega\subset R^3$ with the generalized
vorticity-slip boundary condition (\ref{VSg}). Our main result is the following estimate on the rate of convergence for the solutions:\\
\\
\textbf{Theorem 1.1.} Let $u(t)$ be the smooth solution of the boundary value problem of
Euler equation (\ref{Eu}), (\ref{Eu1}), (\ref{Sl}) on $[0,T]$ with $u(0) = u_0\in H^3(\Omega)$
which satisfies the corresponding assumptions on the initial data in Theorem 1 of \cite{Mas1} (see also Lemma 2.1 in next section), and
$u^\varepsilon(t)$ be the solution of the boundary value problem of the Navier-Stokes equation
(\ref{NS1}), (\ref{NS2}), (\ref{VSg}) with the same initial data $u_0$. Then, there is $T_0>0$ such that
\begin{equation}\label{T1}
  \|u^\varepsilon-u\|^{2}+
  \varepsilon\int_0^{T_1}\|u^\varepsilon-u\|_1^{2}dt
 \leq c\varepsilon^{2-s},\ {\rm on}\ [0,T_0]
 \end{equation}
\begin{equation}\label{T2}
  \|u^\varepsilon-u\|_1^{2} +
  \varepsilon\int_0^{T_1}\|u^\varepsilon-u\|_2^{2}dt
 \leq c\varepsilon^{1-s},\ {\rm on}\ [0,T_0]
 \end{equation}
for any $s>0$ and $\varepsilon$ small enough. Consequently,
\begin{equation}\label{T3}
  \|u^\varepsilon-u\|_{1,p}^{p}\leq c_p\varepsilon^{1-s},\ {\rm on}\ [0,T_0]
 \end{equation}
for $2\leq p<\infty$, any $s>0$ and $\varepsilon$ small enough, and
\begin{equation}\label{T4}
  \|u^\varepsilon-u\|_{L^\infty([0,T]\times\Omega)}
 \leq c\varepsilon^{\frac{2}{5}(1-s')}
 \end{equation}
for any $s' > 0$.

The estimate (\ref{T2}) yields the desired estimates on the rate of  convergence
of the solutions of the Navier-Stokes equations to that of the Euler equations in
$C([0,T], H^1(\Omega))$ norm for general domains with general
Navier-slip boundary conditions. It should be noted that
the estimates (\ref{T1}) and (\ref{T4}) improve (\ref{ke1}) and (\ref{ke2}), and (\ref{T2}), (\ref{T3}) are even better than
(\ref{wxz}) stated in \cite{WXZ} for the special case $B= 0$. Indeed, The estimate
(\ref{T2}) is optimal in the sense that $s$ can not be taken to be $0$ due to boundary layers in general, see Remark 3.1.
below for details.

These estimates are motivated by our recent work \cite{XX3}, where the first order derivative estimates for the solutions are obtained for the Navier-Stokes equations with a
new vorticity boundary condition (1.16), and the work of Masmoudi and Rousset in \cite{Mas1}, where uniform estimates on $\|u^\varepsilon\|_{W^{1.\infty}}$ of the
solutions of the Navier-Stokes equation with the general Navier-slip boundary condition are obtained.

Our approach is an elementary energy estimate for the difference of the solutions between the Navier-Stokes
equations and the Euler equations. Some suitable integrating by part formulas in terms of the vorticity are successfully used to get the optimal rate of convergence.
In comparison to the asymptotic analysis method associate to the boundary layers
(see \cite{IfS,Ke1,Ke2,WXZ}), we should not need any correctors near the boundary.
Meanwhile, the uniform regularity obtained
in \cite{Mas1} plays an essential role in our analysis.

The rest of the paper is organized as follows: In the next section, we give a preliminary
on the well-posedness of the initial boundary value problem of the Navier-Stokes equations (\ref{NS1}),(\ref{NS2})
and (\ref{VSg}), and the local uniform regularity of the solutions.
Then, the estimates on the asymptotic convergence (\ref{T1}), (\ref{T2}), (\ref{T3}) and (\ref{T4}) are established in section 3.

\vskip 1cm

\section{Preliminaries}

Let $\Omega\subset R^3$ be a general bounded smooth domain.
We begin by considering the following Stokes problem:
\begin{eqnarray}\label{St1}
&& \alpha u - \Delta u + \nabla q = f
\textrm{ in }\Omega\\
&& \nabla\cdot u  = 0 \textrm{ in }\Omega\\\label{St2}
&& u\cdot n = 0, n\times(\nabla\times u) = [Bu]_\tau\ {\rm on} \
\partial \Omega
\end{eqnarray}
where $B$ is smooth symmetric tensor.
Set
\[H = \{u\in L^2(\Omega); \nabla \cdot u = 0,  {\rm in} \ \Omega; \ u\cdot n = 0; \ {\rm on}\ \partial\Omega\}\]
\[V = H^1(\Omega)\cap H \]
\[ W  = \{u\in H^2(\Omega); \ n\times (\nabla\times u) = [Bu]_\tau \ {\rm on},
\partial\Omega\}\]
It is well known that for any $v\in V$, one has

\begin{eqnarray}\label{dot}
   \|v\|_1 \leq c \|\nabla\times v\|
\end{eqnarray}

Note that
\[   [Bu]_\tau \cdot \phi = B u\cdot \phi \]
for any $\phi$ satisfying $\phi\cdot n = 0$ on the boundary. We
associate the Stokes problem (\ref{St1})-(\ref{St2}) with the following bilinear form
\begin{eqnarray}
a_\alpha(u,\phi) = \alpha(u, \phi) + \int_{\partial\Omega}B u\cdot \phi + \int_\Omega
(\nabla\times u)\cdot (\nabla\times \phi)
\end{eqnarray}
Note that
\begin{eqnarray}\label{b}
 |\int_{\partial\Omega}B u\cdot u|\leq c\|u\|\|u\|_1
 \end{eqnarray}
for all $u\in V$ and (\ref{dot}).
It follows that $a_\alpha(u,\phi)$ is a positive definite symmetric bilinear form if $\alpha$ large
enough, and is closed with the domain $\mathcal{D}(a_\alpha) = V$. Similar to the discussions
in \cite{XX1}, one can define the self-adjoint operator associated with $a_\alpha$ as
\[  A = \alpha I - P\Delta \]
with the domain $\mathcal{D}(A) = V\cap W$, which implies
\begin{eqnarray}\label{H2}
 \|v\|_2\leq c\|v\| + \|P\Delta v\|, \ \forall v\in \mathcal{D}(A)
 \end{eqnarray}

Now, we turn to the boundary value problem of the Navier-Stokes equations
(\ref{NS1}), (\ref{NS2}) on $\Omega\subset R^3$ with the generalized
vorticity-slip boundary condition (\ref{VSg}).

By using the Galerkin
method based on the orthogonal eigenvectors of $A$, noting
the energy equation
\begin{equation}\label{E1}
\frac{1}{2}\frac{d}{dt}\|u^\varepsilon\|^{2} + \varepsilon
\|\nabla\times u^\varepsilon\|^2 + \varepsilon\int_{\partial\Omega}Bu^\varepsilon\cdot u^\varepsilon
 = 0
\end{equation}
valid for approximate solutions, and (\ref{dot}), (\ref{b}), one can obtain the global existence of weak solutions
to the the boundary value problem of Navier-Stokes equations
(\ref{NS1}), (\ref{NS2}), (\ref{VSg}) (see the corresponding definition in, for instance, \cite{XX1}).

Noting the energy equation
\begin{equation}\label{E2}
\frac{1}{2}\frac{d}{dt}(\|\nabla\times u^\varepsilon\|^{2}+\int_{\partial\Omega}Bu^\varepsilon\cdot u^\varepsilon)  + \varepsilon
\|P\Delta u^\varepsilon\|^2 + (u^\varepsilon\cdot\nabla u^\varepsilon, P\Delta u^\varepsilon)
 = 0
\end{equation}
for approximate solutions, (\ref{dot}), and (\ref{b})-(\ref{E1}), one can show that if $u^\varepsilon_0\in V$,
then there is a maximum time interval $[0,T^\varepsilon)$ such that the weak solutions are
the unique strong one on $[0,T^\varepsilon)$ (see the corresponding definition in, for instance, \cite{XX1}).

It follows from Lemma 3.10 in \cite{XX1} that
\begin{equation}
   (2(S(u)n) - \omega\times n)_\tau = GD(u)_\tau
    \end{equation}
with $GD(u) = -2S(n)u$.

Then, the boundary conditions (\ref{NaS}), (\ref{Nag}) and (\ref{VS}) can be
written to (\ref{VSg}), and (\ref{Nag}) is equivalent to (\ref{VSg}).
Hence, one has the following uniform regularity result for the solutions
to the Navier-Stokes equations:\\
\\
\textbf{Lemma 2.1.}(Masmoudi and Rousset \cite{Mas1}) Let $m>6$ be an integer and
$\Omega$ be a $C^{m+2}$ domain. Consider $u_0\in E^m\cap H$ such that
$\nabla u_0\in W^{1,\infty}_{co}$. Then there is $T_m > 0$ such that for all
sufficient small $\varepsilon$, there is a unique solution $u^\varepsilon\in C([0,T_m],E^m)$
to the Navier-Stokes problem (\ref{NS1}), (\ref{NS2}), (\ref{VSg})
with $u^\varepsilon(0) = u_0$. Moreover, there is a constant $C$ such that
\begin{equation}\label{Ma}
  \|u^\varepsilon\|_{H_{co}^{m}(\Omega)} + \|\nabla u^\varepsilon\|_{H_{co}^{m-1}(\Omega)} +
 \|\nabla u^\varepsilon\|_{W_{co}^{1,\infty}(\Omega)} +
 \varepsilon\int_0^t\|\nabla^2 u^\varepsilon\|_{H_{co}^{m-1}(\Omega)}^2dt\leq C
 \end{equation}
on $[0,T_m]$.\\
\\
Where $H_{co}^{m}(\Omega)$ and $W_{co}^{1,\infty}(\Omega)$ are co-norm vector-spaces,
$\|\cdot\|_{H_{co}^{m}(\Omega)}$ and $\|\cdot\|_{W_{co}^{1,\infty}(\Omega)}$ denote the
corresponding norms, and
\[   E^m = \{u\in H_{co}^{m}(\Omega)| \nabla u \in H_{co}^{m-1}(\Omega)\} \]
Since the notations are rather complicate to be expressed, we omit it here and refer to
Masmoudi and Rousset \cite{Mas1} for the details.
This uniform regularity implies in particular the following uniform bound
\[  \|u^\varepsilon\|_{W^{1,\infty}}\leq C, \ {\rm on} \ [0,T_m], \]
which plays an essential role in our estimates.

\vskip 1cm

\section{Convergence of the solutions}

\vskip 0.5cm

We now turn to the purpose of this paper to establish the convergence
with a rate for the solutions $u^\varepsilon$ to $u$. We start with
the basic $L^2-$estimate.\\
\\
\textbf{Theorem 3.1.} Let $u_{0}\in H^3(\Omega)$ satisfy the assumptions stated in Theorem 1.1.
and $u(t)$ be the solution to the Euler equations (\ref{Eu}), (\ref{Eu1}), (\ref{Sl})
on $[0,T]$ with $u(0) = u_0$, and
$u(t)=u^\varepsilon(t)$ be the solution to the Navier-Stokes problem (\ref{NS1}), (\ref{NS2}), (\ref{VSg})
with $u^\varepsilon(0) = u_0$. Then
\begin{equation}\label{e0}
  \|u^\varepsilon-u\|^{2} +
  \varepsilon\int_0^{T_0}\|u^\varepsilon-u\|_1^{2}dt
 \leq c\varepsilon^{2-s}\ {\rm on}\ [0,T_0]
    \end{equation}
for any $s>0$ and $\varepsilon$ small enough, where $T_0 = \min\{T,T_m\}$.
Consequently,
\begin{equation}
  \|u^\varepsilon-u\|_{L^\infty([0,T]\times\Omega)}
 \leq c\varepsilon^{\frac{2}{5}(1-s')}
 \end{equation}
for any $s' > 0$.\\
\\
\textbf{Proof:} Note that $u^\varepsilon-u$ satisfies
\begin{eqnarray}\label{D1}
&& \partial_{t}(u^\varepsilon-u) -\varepsilon \Delta (u^\varepsilon-u)
   + \Phi  + \nabla (p^\varepsilon-p) = \varepsilon\Delta u,
\textrm{ in }\Omega\\ \label{D2}
&& \nabla\cdot u^\varepsilon  = 0, \textrm{ in }\Omega\\\label{D3}
&&(u^\varepsilon-u)\cdot n = 0,\  n\times(\omega^\varepsilon-\omega) = [B(u^\varepsilon-u)+B u]_\tau -n\times\omega  \ {\rm on} \
\partial \Omega
\end{eqnarray}
where $\omega= \nabla\times u$, $\omega^\varepsilon= \nabla\times u^\varepsilon$,
\[ \Phi = u\cdot\nabla(u^\varepsilon-u) + (u^\varepsilon-u)\cdot\nabla u +(u^\varepsilon-u)\cdot\nabla(u^\varepsilon-u)
  \]
Then, the following identity holds
\begin{equation}\label{de1}
\frac{1}{2}\frac{d}{dt}\|u^\varepsilon-u\|^{2} + \varepsilon
\|\nabla\times(u^\varepsilon-u)\|^2 + \mathcal{B}_0
  + (\Phi,u^\varepsilon-u) = \varepsilon(\Delta u,u^\varepsilon-u)
\end{equation}
where
\[ \mathcal{B}_0 = \varepsilon\int_{\partial\Omega}n\times(\omega^\varepsilon-\omega)(u^\varepsilon-u)
= \varepsilon\int_{\partial\Omega}(B(u^\varepsilon-u)+ Bu-n\times\omega)(u^\varepsilon-u) \]
Note that
\[ \int_{\partial\Omega}n\times(\omega^\varepsilon-\omega)(u^\varepsilon-u)
=\int_{\partial\Omega}(B(u^\varepsilon-u)+ Bu-n\times\omega)(u^\varepsilon-u)\]
\[  \leq  c\int_{\partial\Omega}(|(u^\varepsilon-u)|^2+|(u^\varepsilon-u)|) \]
\[ \leq c(\|u^\varepsilon-u\|\|\omega^\varepsilon-\omega\| + |u^\varepsilon-u|_{L(\partial\Omega)})\]
It follows from the trace theorem that
\begin{equation}\label{Tr}
|u^\varepsilon-u|_{L(\partial\Omega)}\leq |u^\varepsilon-u|_{L^q(\partial\Omega)}
 \leq c\|u^\varepsilon-u\|_{H^s(\Omega)}
\end{equation}
for any $q>1$ and then any $s>0$.

By interpolation (see \cite{Tem}), we have
\begin{equation}\label{Intp}
\|u^\varepsilon-u\|_{H^s(\Omega)}\leq c\|u^\varepsilon-u\|^{(1-s)}\|u^\varepsilon-u\|_1^{s}
 \leq c\|u^\varepsilon-u\|^{(1-s)}\|\omega^\varepsilon-\omega\|^s
\end{equation}
Note that
\begin{equation}\label{2.1}
\varepsilon\|u^\varepsilon-u\|^{(1-s)}\|\omega^\varepsilon-\omega\|^s
 \leq \delta\varepsilon\|\omega^\varepsilon-\omega\|^2 + c_\delta(\|u^\varepsilon-u\|^2 + \varepsilon^{2-s})
\end{equation}
for any $s\in (0,1)$ and
\[ \varepsilon \|u^\varepsilon-u\|\|\omega^\varepsilon-\omega\|
\leq \varepsilon^2\|\omega^\varepsilon-\omega\|^2 + \|u^\varepsilon-u\|^2 \]
Then, it holds that
\[ \mathcal{B}_0\leq
 2\delta\varepsilon \|\omega^\varepsilon-\omega\|^2 + c_\delta\|u^\varepsilon-u\|^2 + \varepsilon^{2-s}\]
for any $s\in (0,1)$ and $\varepsilon$ small enough.\\

Note also that
\[ (\Phi,u^\varepsilon-u) = ((u^\varepsilon-u)\cdot\nabla u, u^\varepsilon-u)\leq  c\|u^\varepsilon-u\|^{2}, \]
and
\[ \varepsilon(\Delta u,u^\varepsilon-u) \leq \|u^\varepsilon-u\|^{2} + c\varepsilon^2 \]
It follows that
\begin{equation}\label{ee}
\frac{d}{dt}\|u^\varepsilon-u\|^{2} + \varepsilon
\|\nabla\times(u^\varepsilon-u)\|^2 \leq
c(\|u^\varepsilon-u\|^2 + \varepsilon^{2-s})
\end{equation}
for any $s\in (0,1)$, and then for any $s>0$ and $\varepsilon$ small enough. Note that
$u^\varepsilon(0)-u(0)=0$.
Then, (\ref{e0}) follows from (\ref{ee}) and the Gronwall lemma.\\
\\
Consequently, by using the Gagliardo-Nirenberg interpolation inequality, one has
\begin{equation}
  \|u^\varepsilon-u\|_{L^\infty(\Omega)}
  \leq c\|u^\varepsilon-u\|^{\frac{2}{5}}\|u^\varepsilon-u\|_{W^{1,\infty}(\Omega)}^{\frac{3}{5}}
 \leq c\varepsilon^{\frac{2}{5}(1-s')}
 \end{equation}
for any $s' > 0$. The theorem is proved.\\
\\
Next, we prove the major estimate in this paper.\\
\\
\textbf{Theorem 3.2.} Let $u_{0}\in H^3(\Omega)$ satisfy the assumptions stated in Theorem 1.1.,
and $u(t)$ be the solution to the Euler equations (\ref{Eu}), (\ref{Eu1}), (\ref{Sl})
on $[0,T]$ with $u(0) = u_0$, and
$u(t)=u^\varepsilon(t)$ be the solution to the Navier-Stokes problem (\ref{NS1}), (\ref{NS2}), (\ref{VSg})
with $u^\varepsilon(0) = u_0$. Then
\begin{equation}
  \|u^\varepsilon-u\|_1^{2} +
  \varepsilon\int_0^T\|u^\varepsilon-u\|_2^{2}dt
 \leq c\varepsilon^{1-s} \ {\rm on}\ [0,T_0]
    \end{equation}
for any $s>0$ and $\varepsilon$ small enough, where $T_0 = \min\{T,T_m\}$.\\
\\
\textbf{Proof:} Let $s$ be the same as in Theorem 3.1..
Note the smoothness of the solutions and
\[   \partial_t (u^\varepsilon-u)\cdot n = 0\ \rm{on} \ \partial\Omega \]
It follows from (\ref{D1}),(\ref{D2}) and (\ref{D3}) that
\[
\frac{1}{2}\frac{d}{dt}\|(\omega^\varepsilon-\omega)\|^{2} + \varepsilon
\|P\Delta(u^\varepsilon-u)\|^2 -(\Phi,
  P\Delta(u^\varepsilon-u))\]
\[ = \int_{\partial\Omega} \partial_t (u^\varepsilon-u)\cdot(n\times(\omega^\varepsilon-\omega))
 - \varepsilon(\Delta u,P\Delta(u^\varepsilon-u))
\]
where $P$ is the Lerray projection,
\[ \Phi = u\cdot\nabla(u^\varepsilon-u)+ (u^\varepsilon-u)\cdot\nabla u + (u^\varepsilon-u)\cdot\nabla (u^\varepsilon-u) \]
and that
\[ \int_{\partial\Omega} \partial_t (u^\varepsilon-u)\cdot(n\times(\omega^\varepsilon-\omega))=
\int_{\partial\Omega} \partial_t (u^\varepsilon-u)\cdot(B (u^\varepsilon-u) + Bu -  n\times\omega) \]
\[ = \frac{1}{2}\frac{d}{dt}(\int_{\partial\Omega}B(u^\varepsilon-u)\cdot (u^\varepsilon-u)  +
 2\int_{\partial\Omega}  (u^\varepsilon-u)\cdot(Bu -  n\times\omega)) - \mathcal{B}_1\]
where
\[ \mathcal{B}_1
   = \int_{\partial\Omega}(u^\varepsilon-u)\cdot\partial_t(Bu -  n\times\omega)\]
It follows from Theorem 3.1. and (\ref{2.1})
\[ |\mathcal{B}_1|\leq c\int_{\partial\Omega}|(u^\varepsilon-u)|
\leq \delta\|(\omega^\varepsilon-\omega)\|^{2} +  c\varepsilon^{1-s}\]
Note also that
\[  \varepsilon|(\Delta u,P\Delta(u^\varepsilon-u))|\leq \frac{\varepsilon}{2}\|P\Delta(u^\varepsilon-u)\|^2+ c\varepsilon \]
and that
\[ - (\Phi,P\Delta(u^\varepsilon-u))= (P\Phi,-\Delta(u^\varepsilon-u))  \]
\[  = (\nabla\times\Phi,(\omega^\varepsilon-\omega))
+ \int_{\partial\Omega}n\times(\omega^\varepsilon-\omega)\cdot P\Phi \]
\[  = (\nabla\times\Phi,(\omega^\varepsilon-\omega))
+ \int_{\partial\Omega}(B(u^\varepsilon-u) + Bu- n\times\omega)\cdot P\Phi \]
Then, we get
\begin{equation}\label{EI}
  \frac{1}{2}\frac{d}{dt} E
  + \frac{\varepsilon}{2}\|P\Delta(u^\varepsilon-u)\|^2 \leq N + BN + BNL
 + c(\|(\omega^\varepsilon-\omega)\|^{2} +  \varepsilon^{1-s})
    \end{equation}
for $\varepsilon$ small enough, where
\[ E = \|(\omega^\varepsilon-\omega)\|^{2}-(\int_{\partial\Omega}B(u^\varepsilon-u)\cdot (u^\varepsilon-u)  +
 2\int_{\partial\Omega}  (u^\varepsilon-u)\cdot(Bu -  n\times\omega)) \]
\[  N = |(\nabla\times\Phi,(\omega^\varepsilon-\omega))|  \]
\[ BN = |\int_{\partial\Omega}B(u^\varepsilon-u)\cdot P\Phi| \]
\[ BNL = |\int_{\partial\Omega}( Bu- n\times\omega)\cdot P\Phi| \]

The term $N$ can be estimated easily by using the Sobolev inequalities
and the known uniform bounds for $\|u^\varepsilon\|_\infty$ and $\|\nabla u^\varepsilon\|_\infty$.
The term $BN$ can be estimated after integrating by part properly. While, the estimate for the leading order
term $BNL$ is  rather complicated due to the possible appearance of boundary layers. We now
carry out these estimates.\\
\\
\textbf{Estimates on $N$:}\\
\\
The term $N$ is estimated as follows:\\
\\
Note that
\[  \nabla\times\Phi = \nabla\times(\omega\times(u^\varepsilon-u) + (\omega^\varepsilon-\omega)\times u
+ (\omega^\varepsilon-\omega)\times (u^\varepsilon-u) \]
\[ = [\omega, u^\varepsilon-u] + [\omega^\varepsilon-\omega, u] + [\omega^\varepsilon-\omega,\omega^\varepsilon-\omega] \]
where
\[ [\varphi,\psi] = \psi\cdot\nabla\varphi - \varphi\cdot\nabla \psi \]
And
\[ |((u^\varepsilon-u)\cdot\nabla \omega,\omega^\varepsilon-\omega)|\leq c\|\omega^\varepsilon-\omega\|^2 \]
\[ |(\omega\cdot\nabla(u^\varepsilon-u),\omega^\varepsilon-\omega)|
\leq c\|u^\varepsilon-u\|_1\|\omega^\varepsilon-\omega\|\leq c\|\omega^\varepsilon-\omega\|^2 \]
It follows that
\[ |([\omega, u^\varepsilon-u],\omega^\varepsilon-\omega)| \leq c\|\omega^\varepsilon-\omega\|^2 \]
On the other hand,
\[ (u\cdot\nabla(\omega^\varepsilon-\omega),\omega^\varepsilon-\omega) = 0 \]
\[ ((\omega^\varepsilon-\omega)\nabla u,\omega^\varepsilon-\omega)\leq c\|\omega^\varepsilon-\omega\|^2 \]
It follows that
\[ |([\omega^\varepsilon-\omega, u],\omega^\varepsilon-\omega)| \leq c\|\omega^\varepsilon-\omega\|^2 \]
Note that
\[ ((u^\varepsilon-u)\cdot\nabla(\omega^\varepsilon-\omega),\omega^\varepsilon-\omega) = 0 \]
\[ ((\omega^\varepsilon-\omega)\cdot\nabla(u^\varepsilon-u),\omega^\varepsilon-\omega) \]
\[ \leq \|\nabla(u^\varepsilon-u)\|_\infty\|\omega^\varepsilon-\omega\|^2\leq c\|\omega^\varepsilon-\omega\|^2 \]
(here and below we will use the uniform regularity that $\|\nabla u^\varepsilon\|_\infty\leq c$).
It follows that
\[ |([\omega^\varepsilon-\omega, u^\varepsilon-u],\omega^\varepsilon-\omega)| \leq c\|\omega^\varepsilon-\omega\|^2 \]
Hence
\begin{equation}\label{N}
   N \leq c\|\omega^\varepsilon-\omega\|^2
    \end{equation}
\\
\textbf{Estimates on $BN$:}\\

Next, we estimate the term
\[  BN = |\int_{\partial\Omega}B(u^\varepsilon-u)\cdot P\Phi| \]
Note that
\[  P\Phi = \Phi + \nabla p_\Phi \]
involves a scalar function $p_\Phi$, which is difficult to estimate on the boundary.
We will transform it to an estimate on $\Omega$
by an integrating by part formula.

To this end, we first extend $n, B$ to the interior of $\Omega$ as follows:
\[  n(x) = \varphi(r(x))\nabla(r(x)) \]
\[  B(x) = \varphi (r(x))B(\Pi x) \]
where
\[ r(x) = \min_{y\in\partial\Omega}d(x,y) \]
\[ \Pi x = y_x\in \partial\Omega \]
such that
\[  r(x) = d(x,y_x) \]
which is uniquely defined on $\Omega_\sigma = \{x\in \Omega,\ r(x)\leq 2\sigma\}$ for some $\sigma>0$,
and $\varphi(s)$ is smooth and compactly supported in $[0, 2\sigma)$
such that
\[  \varphi(0) = 1, \ {\rm on} \ [0,\sigma] \]

Then, we can deduce the estimate of the boundary term $BN$ to an interior estimate on $\Omega$
by the Stokes formula
\[ \int_{\partial\Omega}B(u^\varepsilon-u)\cdot P\Phi
 = \int_{\partial\Omega}(n\times B(u^\varepsilon-u)\cdot(n\times P\Phi) \]
\[ = (n\times B(u^\varepsilon-u), \nabla\times\Phi)
-(\nabla\times(n\times B(u^\varepsilon-u)),P\Phi)
  \]
since $P\Phi\cdot n = 0$ on the boundary, and $\nabla\times P\Phi = \nabla\times\Phi$.

Note that
\[ |(n\times B(u^\varepsilon-u), (u^\varepsilon-u)\cdot\nabla \omega)|
\leq c\|u^\varepsilon-u\|^2 \leq c\|\omega^\varepsilon-\omega\|^2 \]
\[ |(n\times B(u^\varepsilon-u), \omega\cdot\nabla(u^\varepsilon-u))|
\leq c\|u^\varepsilon-u\|_1^2 \leq c\|\omega^\varepsilon-\omega\|^2  \]
It follows that
\[ |(n\times B(u^\varepsilon-u),[\omega, u^\varepsilon-u])| \leq c\|\omega^\varepsilon-\omega\|^2  \]
Note that $u\cdot n = 0$ on $\partial\Omega$, $\nabla\cdot u = 0$. Then
\[ |(n\times B(u^\varepsilon-u), u\cdot\nabla(\omega^\varepsilon-\omega))|
\]
\[ = |(\omega^\varepsilon-\omega, u\cdot\nabla(n\times B(u^\varepsilon-u)))|
\leq c\|u^\varepsilon-u\|_1^2\leq c\|\omega^\varepsilon-\omega\|^2 \]
and
\[ |(n\times B(u^\varepsilon-u),(\omega^\varepsilon-\omega)\cdot\nabla u)|\leq c\|\omega^\varepsilon-\omega\|^2 \]
It follows that
\[ |(n\times B(u^\varepsilon-u),[\omega^\varepsilon-\omega,u])| \leq c\|\omega^\varepsilon-\omega\|^2  \]
Similarly
\[ |(n\times B(u^\varepsilon-u), (u^\varepsilon-u)\cdot\nabla(\omega^\varepsilon-\omega))|
= |(\omega^\varepsilon-\omega, (u^\varepsilon-u)\cdot\nabla(n\times B(u^\varepsilon-u)))| \]
\[\leq c\|u^\varepsilon-u\|_1^2\leq c\|\omega^\varepsilon-\omega\|^2 \]
and
\[ |(n\times B(u^\varepsilon-u),(\omega^\varepsilon-\omega)\cdot\nabla (u^\varepsilon-u))|\leq c\|\omega^\varepsilon-\omega\|^2 \]
It follows that
\[ |(n\times B(u^\varepsilon-u),[\omega^\varepsilon-\omega,u^\varepsilon-u])| \leq c\|\omega^\varepsilon-\omega\|^2  \]
Hence
\begin{equation}\label{II}
  |(n\times B(u^\varepsilon-u), \nabla\times\Phi)| \leq c\|\omega^\varepsilon-\omega\|^2
    \end{equation}
\\
Note also that
\[ |(\nabla\times(n\times B(u^\varepsilon-u)), P\Phi)|
\leq \|\nabla\times(n\times B(u^\varepsilon-u))\|\|\Phi\| \]
\[ \leq c\|(u^\varepsilon-u))\|_1^2
\leq c\|\omega^\varepsilon-\omega\|^2 \]
Hence, we have
\begin{equation}\label{BN}
   BN \leq c\|\omega^\varepsilon-\omega\|^2
    \end{equation}
\\
\textbf{Estimate on $\int_{\partial\Omega}( Bu- n\times\omega)\cdot P\Phi$:}\\
\\
Finally, we estimate the leading order term on the boundary
\[ BNL = |\int_{\partial\Omega}( Bu- n\times\omega)\cdot P\Phi| \]
It should be noted that the estimate is trivial if
\[   [Bu]_\tau- n\times\omega = 0 \]
which is that the Euler solution satisfies the same boundary condition (1.11) just as the Navier-Stokes solution does so there is no strong boundary layer. This applies to the cases treated in \cite{XX,XX2}.
Here, $[Bu]_\tau- n\times\omega$ may be not equal zero, then
boundary layers may occur. Additional efforts are needed to overcome the new difficulties.\\
\\
Similar to the above, we have
\[ |\int_{\partial\Omega}( Bu- n\times\omega)\cdot P\Phi| =
|\int_{\partial\Omega}n\times( Bu- n\times\omega)\cdot n\times P\Phi| \]
\[ = |(n\times ( Bu- n\times\omega), \nabla\times\Phi)
-(\nabla\times(n\times ( Bu- n\times\omega)),P\Phi)|
  \]
The term $|(n\times ( Bu- n\times\omega), \nabla\times\Phi)|$ is relatively easy to
be estimated since $\nabla\times\Phi$ does not involve the pressure function.\\
\\
Rewrite $\Phi$ as
\[  \Phi = \Phi_1 + \Phi_2 + \Phi_3  \]
where
\[ \Phi_1 = (u^\varepsilon-u)\cdot\nabla u  \]
\[ \Phi_2 = u\cdot\nabla(u^\varepsilon-u) \]
\[ \Phi_3 = (u^\varepsilon-u)\cdot\nabla(u^\varepsilon-u) \]
Note that
\[ |(n\times ( Bu- n\times\omega), \nabla\times\Phi_1)| = |(n\times ( Bu- n\times\omega), \nabla\times((u^\varepsilon-u)\cdot\nabla u))| \]
\[ = |\int_{\partial\Omega} n\times((u^\varepsilon-u)\cdot\nabla u)\cdot (n\times ( Bu- n\times\omega))
+ (\nabla\times(n\times ( Bu- n\times\omega)),(u^\varepsilon-u)\cdot\nabla u)| \]
\[\leq c(\int_{\partial\Omega}|u^\varepsilon-u| + \|u^\varepsilon-u\|) \leq c(\|\omega^\varepsilon-\omega\|^2 + \varepsilon^{1-s}) \]
and
\[ |(n\times ( Bu- n\times\omega), \nabla\times\Phi_3)| =
|(n\times ( Bu- n\times\omega), \nabla\times((u^\varepsilon-u)\cdot\nabla (u^\varepsilon-u)))| \]
\[ = |\int_{\partial\Omega} n\times((u^\varepsilon-u)\cdot\nabla (u^\varepsilon-u))\cdot (n\times ( Bu- n\times\omega))
+ (\nabla\times(n\times ( Bu- n\times\omega)),(u^\varepsilon-u)\cdot\nabla (u^\varepsilon-u))| \]
\[\leq c(\int_{\partial\Omega}|u^\varepsilon-u| + \|u^\varepsilon-u\|_1^2) \leq c(\|\omega^\varepsilon-\omega\|^2 + \varepsilon^{1-s}) \]
since $u^\varepsilon-u$ is smooth, so that
\[  |\nabla(u^\varepsilon-u)| \leq c \|\nabla(u^\varepsilon-u)\|_\infty \leq c, \ {\rm on}\ \partial\Omega \]
To estimate $(n\times ( Bu- n\times\omega), \nabla\times\Phi_2)$, we note
\[ \nabla\times(\psi\cdot\nabla\phi) = \psi\cdot\nabla(\nabla\times\phi) + \nabla\psi^\perp\cdot\nabla\phi\]
where $\nabla\psi^\perp$ is expressed in component by
\begin{equation}\label{ex}
(\nabla\psi^\perp\cdot\nabla\phi)_j = (-1)^{j+1}\partial_{j+1}\psi\cdot\nabla\phi_{j+1} + (-1)^{j+2}\partial_{j+2}\psi\cdot\nabla\phi_{j+2}
\end{equation}
with the index modulated by 3. Hence, we have
\[ |(n\times ( Bu- n\times\omega), \nabla\times\Phi_2)| = |(n\times ( Bu- n\times\omega), \nabla\times(u\cdot\nabla(u^\varepsilon-u)))| \]
\[ = |(n\times ( Bu- n\times\omega), u\cdot\nabla(\omega^\varepsilon-\omega))
+ (n\times ( Bu- n\times\omega), \nabla u^\perp\cdot\nabla(u^\varepsilon-u))| \]
Note that
\[ |(n\times ( Bu- n\times\omega), u\cdot\nabla(\omega^\varepsilon-\omega)) |
= |(\omega^\varepsilon-\omega, u\cdot\nabla(n\times ( Bu- n\times\omega)))| \]
\[ = (\nabla\times(u^\varepsilon-u), u\cdot\nabla(n\times ( Bu- n\times\omega))) \]
\[ = |\int_{\partial\Omega} n\times(u^\varepsilon-u)\cdot(u\cdot\nabla(n\times ( Bu- n\times\omega)))
+ (u^\varepsilon-u, \nabla\times(u\cdot\nabla(n\times ( Bu- n\times\omega))))| \]
\[\leq c(\int_{\partial\Omega}|u^\varepsilon-u| + \|u^\varepsilon-u\|) \leq c(\|\omega^\varepsilon-\omega\|^2 + \varepsilon^{1-s}) \]
Note also that
\[ (\partial_j u\cdot\nabla(u^\varepsilon-u)_j,(n\times ( Bu- n\times\omega))_k)
=  \]
\[(\partial_j u,\nabla((u^\varepsilon-u)_j(n\times ( Bu- n\times\omega))_k)
-(\partial_j u\cdot\nabla(n\times ( Bu- n\times\omega))_k,(u^\varepsilon-u)_j),\]
and that
\[ |(\partial_j u,\nabla((u^\varepsilon-u)_j(n\times ( Bu- n\times\omega))_k)|
= |\int_{\partial\Omega}(u^\varepsilon-u)_j(n\times ( Bu- n\times\omega))_k\partial_j u\cdot n |\]
\[\leq c\|u\|_{1,\infty}^2\int_{\partial\Omega}|u^\varepsilon-u| \leq c\int_{\partial\Omega}|u^\varepsilon-u|
 \leq c(\|\omega^\varepsilon-\omega\|^2 + \varepsilon^{1-s}) \]
since $\nabla\cdot \partial_j u = 0$, and that
\[|(\partial_j u\cdot\nabla(n\times ( Bu- n\times\omega))_k,(u^\varepsilon-u)_j)|
\leq c\|u^\varepsilon-u\|\leq c\varepsilon^{1-s} \]
it follows that
\[ |(n\times ( Bu- n\times\omega), \nabla\times\Phi_2)|\leq c(\|\omega^\varepsilon-\omega\|^2 + \varepsilon^{1-s}) \]
Then, we get
\begin{equation}
 |(n\times ( Bu- n\times\omega), \nabla\times\Phi)| \leq c(\|\omega^\varepsilon-\omega\|^2 + \varepsilon^{1-s})
 \end{equation}
\\
It remains to estimate $|(\nabla\times(n\times(Bu- n\times\omega)),P\Phi)|$. There is also a
difficulty arising from $P\Phi$. To do it, it is noticed that
\[ P\Phi = u\cdot\nabla(u^\varepsilon-u)- (u^\varepsilon-u)\cdot\nabla u + P\tilde{\Phi} \]
where
\[  \tilde{\Phi} = 2(u^\varepsilon-u)\cdot\nabla u + (u^\varepsilon-u)\cdot\nabla (u^\varepsilon-u) \]
with
\[   u\cdot\nabla(u^\varepsilon-u)- (u^\varepsilon-u)\cdot\nabla u = \nabla\times((u^\varepsilon-u)\times u) \in H \]
since that $(u^\varepsilon-u)\cdot n = 0$ and $ u\cdot n = 0$ on $\partial\Omega$ will imply
\[  (u^\varepsilon-u)\times u = \lambda n \]
and
\[  H = \{v = \nabla\times \phi;\ \phi\in H^1(\Omega), n\times \phi = 0\} \]
Hence, we have
\[   P(u\cdot\nabla(u^\varepsilon-u)- (u^\varepsilon-u)\cdot\nabla u) = u\cdot\nabla(u^\varepsilon-u)- (u^\varepsilon-u)\cdot\nabla u \]
so that
\[ (\nabla\times(n\times(Bu- n\times\omega)),P\Phi)=  \]
\[ (\nabla\times(n\times(Bu- n\times\omega)),u\cdot\nabla(u^\varepsilon-u)- (u^\varepsilon-u)\cdot\nabla u)
 + (\nabla\times(n\times(Bu- n\times\omega)),P\tilde{\Phi}) \]
and
\[  \|P\tilde{\Phi}\| \leq c\|\tilde{\Phi}\|
\leq c(\|u\|_{1,\infty}+\|u^\varepsilon-u\|_{1,\infty})\|u^\varepsilon-u\|
\leq c \|u^\varepsilon-u\| \leq c\varepsilon^{1-s}\]
These properties help us to complete the rest of estimates.

First, we have
\[ |(\nabla\times(n\times(Bu- n\times\omega)),P\tilde{\Phi})|
\leq c\|P\tilde{\Phi}\| \leq c\|\tilde{\Phi}\| \leq c \|u^\varepsilon-u\|\leq c\varepsilon^{1-s}  \]
Next, note that
\[ |(\nabla\times(\times(Bu- n\times\omega)),u\cdot\nabla(u^\varepsilon-u))
 = \]
\[ |(u^\varepsilon-u,u\cdot\nabla(\nabla\times(\times(Bu- n\times\omega))))|
 \leq c\|u\|_3^2\|u^\varepsilon-u\|\leq c\varepsilon^{1-s} \]
since $\nabla u = 0$ in $\Omega$ and $u\cdot n = 0$ on $\partial\Omega$.

Note also that
\[|(\nabla\times(n\times(Bu- n\times\omega)) , (u^\varepsilon-u)\nabla u)|
\leq c\|u\|_2^2\|u^\varepsilon-u\|\leq c\varepsilon^{1-s}  \]
Hence, we have
\begin{equation}\label{V}
  |(\nabla\times(n\times(Bu- n\times\omega)),P\Phi)|\leq c\varepsilon^{1-s}
    \end{equation}
Then, we conclude that
\begin{equation}\label{BNL}
   BNL = |\int_{\partial\Omega}( Bu- n\times\omega)\cdot P\Phi|\leq c(\|\omega^\varepsilon-\omega\|^2 + \varepsilon^{1-s})
    \end{equation}
\\
In conclusion, it follows from (\ref{EI}),(\ref{N}),(\ref{BN}) and (\ref{BNL}) that
\begin{equation}
 \frac{1}{2}\frac{d}{dt} E
  + \frac{\varepsilon}{2}\|P\Delta(u^\varepsilon-u)\|^2 \leq
  c(\|\omega^\varepsilon-\omega\|^{2} + \varepsilon^{1-s})
    \end{equation}

Recall that
\[ E = \|(\omega^\varepsilon-\omega)\|^{2}-(\int_{\partial\Omega}B(u^\varepsilon-u)\cdot (u^\varepsilon-u)  +
 2\int_{\partial\Omega}(u^\varepsilon-u)\cdot(Bu -  n\times\omega)) \]
and
\[|\int_{\partial\Omega}B(u^\varepsilon-u)\cdot (u^\varepsilon-u)|\leq c\int_{\partial\Omega}|(u^\varepsilon-u)|^2\leq
c\|u^\varepsilon-u\|\|\omega^\varepsilon-\omega\|
\leq  \|\omega^\varepsilon-\omega\|^2 + c\|u^\varepsilon-u\|^2 \]
and
\[  |\int_{\partial\Omega}(u^\varepsilon-u)\cdot(Bu -  n\times\omega)|
\leq c\int_{\partial\Omega}|(u^\varepsilon-u)|
\leq\|\omega^\varepsilon-\omega\|^2 + c\|u^\varepsilon-u\|^2 \]
It follows that
\begin{equation}
\|\omega^\varepsilon-\omega\|^2
  + \varepsilon c\int_0^t \|P\Delta(u^\varepsilon-u)\|^2 \leq
  c\int_0^t\|\omega^\varepsilon-\omega\|^{2} +  c\varepsilon^{1-s}
    \end{equation}
By using the Gronwall Lemma, one gets that
\begin{equation}
\|\omega^\varepsilon-\omega\|^2\leq c\varepsilon^{1-s}
    \end{equation}
on $[0, T_0]$ for any given $s>0$ and $\varepsilon$ for small enough.
Then
\begin{equation}
\|u^\varepsilon-u\|_1^2\leq c\varepsilon^{1-s}
    \end{equation}
and
\begin{equation}
  \varepsilon \int_0^t \|P\Delta(u^\varepsilon-u)\|^2 \leq
    c\varepsilon^{1-s}
    \end{equation}
Combining this with (\ref{H2}) and Theorem 3.1 implies that
\begin{equation}
  \varepsilon \int_0^t \|(u^\varepsilon-u)\|_2^2 \leq
    c\varepsilon^{1-s}
    \end{equation}
Note that
\[   \|\nabla(u^\varepsilon-u)\|_p^p\leq c\|\nabla(u^\varepsilon-u)\|_\infty^{p-2}\|\nabla(u^\varepsilon-u)\|^2 \]
It follows that
\[   \|\nabla(u^\varepsilon-u)\|_p^p\leq  c\varepsilon^{1-s} \]
The theorem is proved.\\
\\
\textbf{Remark 3.1.} The estimate (\ref{T2}) is optimal in the sense that $s$ can not be taken to be $0$.
Since if $s=0$, then (1.19) will imply that $u^\varepsilon- u$ is uniform bounded in $L^2(0,T; H^2(\Omega))$. Then,
there is a subsequence such that $u^{\varepsilon_n}- u$ weakly convergence in $H^2(\Omega)$ for
a.e. $t\in [0,T]$. Note that $u^{\varepsilon_n}- u$ strongly convergence to $0$ in $H^1(\Omega)$
for all $t\in [0,T]$. It follows that
\[    u^{\varepsilon_n}\rightarrow u\ weakly\ in\ H^2(\Omega),\ a.e. \ t \]
Note that $W$ is closed in $H^2(\Omega)$, then is weakly closed. Hence, $u\in W$ for
a.e. $t\in [0,T]$, i.e.,
\[   [Bu]_\tau- n\times\omega = 0,\ \ a.e. \ t \]
Since $u$ is smooth, then it holds for all $t\in [0,T]$. This is impossible in general, see \cite{BF2,BF3}.

\end{document}